# Ratio Estimators In Simple Random Sampling When Study Variable Is An Attribute

Rajesh Singh  and Mukesh Kumar

Department of Statistics, Banaras Hindu University, Varanasi-221005, INDIA

Florentin Smarandache

Department of Mathematics, University of New Mexico, Gallup, USA

rsinghstat@yahoo.com, mukesh.stat@gmail.com, smarand@unm.edu

**Abstract**

In this paper we have suggested a family of estimators for the population mean when study variable itself is qualitative in nature. Expressions for the bias and mean square error (MSE) of the suggested family have been obtained. An empirical study has been carried out to show the superiority of the constructed estimator over others.

## 1. Introduction

The use of auxiliary information can increase the precision of an estimator when study variable y is highly correlated with auxiliary variable x. In many situations study variable is generally ignored not only by ratio scale variables that are essentially qualitative , or  nominal scale, in nature, such as sex, race, colour, religion, nationality, geographical region, political upheavals (see [1]). Taking into consideration the point biserial correlation coefficient between auxiliary attribute and study variable, several authors including [2], [3], [4], [5] and [6] defined ratio estimators of population mean when the priori information of population proportion of units, possessing some attribute is available. All the others have implicitly assumed that the study variable Y is quantitative whereas the auxiliary variable is qualitative.

In this paper we consider some estimators in which study variable itself is qualitative in nature. For example suppose we want to study the labour force participation (LFP) decision of adult males. Since an adult is either in the labour force or not, LFP is a yes or no decision. Hence, the study variable can take two values, say 1, if the person is in the labour force and 0 if he is not. Labour economics research suggests that the LFP decision is a function of the unemployment rate, average wage rate, education, family income, etc (See [1]).

Consider a sample of size n drawn by simple random sampling without replacement (SRSWOR) from a population size N. Let $\phi_i$ and $x_i$ denote the observations on variable $\phi$ and x respectively for $i^{th}$ unit (i=1,2,3...N). $\phi_i = 1$, if $i^{th}$ unit of population possesses attribute $\phi$ and $\phi_i = 0$, otherwise. Let $A = \sum_{i=1}^{N} \phi_i$ and $a = \sum_{i=1}^{n} \phi_i$ denote the total number of units in the population and sample possessing attribute $\phi$ respectively, $P = \dfrac{A}{N}$ and $p = \dfrac{a}{n}$ denote the proportion of units in the population and sample, respectively, possessing attribute $\phi$.

Define,

$$e_\phi = \frac{(p - P)}{P}, \qquad e_x = \frac{(\bar{x} - \bar{X})}{\bar{X}}$$

Such that,

$$E(e_i) = 0, (i = \phi, x)$$

and

$$E\left(e_\phi^2\right) = fC_p^2, \qquad E\left(e_x^2\right) = fC_x^2, \qquad E\left(e_x e_\phi\right) = f\rho_{pb}C_\phi C_x.$$

where,

$$f = \left(\frac{1}{n} - \frac{1}{N}\right), \qquad C_p^2 = \frac{S_p^2}{P^2}, \qquad C_x^2 = \frac{S_x^2}{\overline{X}^2},$$

and $\rho_{pb} = \dfrac{S_{\phi x}}{S_\phi S_x}$ is the point biserial correlation coefficient.

Here,

$$S_\phi^2 = \frac{1}{N-1}\sum_{i=1}^{N}(\phi_i - P)^2, \qquad S_x^2 = \frac{1}{N-1}\sum_{i=1}^{N}(x_i - \overline{X})^2 \quad \text{and} \quad S_{\phi x} = \frac{1}{N-1}\left(\sum_{i=1}^{N}\phi_i x_i - NP\overline{X}\right).$$

## 2. The proposed estimator

We first propose the following ratio-type estimator

$$t_1 = \left(\frac{P}{\overline{x}}\right)\overline{X} \qquad (2.1)$$

The bias and MSE of the estimator $t_1$, to the first order of approximation is respectively,

given by

$$B(t_1) = f\left(\frac{C_x^2}{2} - \rho_{pb}C_\phi C_x\right) \qquad (2.2)$$

$$MSE(t_1) = f\left(C_\phi^2 + C_x^2 - 2\rho_{pb}C_\phi C_x\right) \qquad (2.3)$$

Following [7], we propose a general family of estimators for P as

$$t_2 = H(p, u)$$ (2.4)

where $u = \dfrac{\overline{x}}{\overline{X}}$ and $H(p, u)$ is a parametric equation of $p$ and $u$ such that

$$H(p, l) = P, \forall P$$ (2.5)

and satisfying following regulations:

(i) Whatever be the sample chosen, the point $(p, u)$ assume values in a bounded closed convex subset $R_2$ of the two-dimensional real space containing the point $(p, l)$.

(ii) The function $H(p, u)$ is a continuous and bounded in $R_2$.

(iii) The first and second order partial derivatives of $H(p, u)$ exist and are continuous as well as bounded in $R_2$.

Expanding $H(p, u)$ about the point $(P, 1)$ in a second order Taylor series we have

$$
\begin{aligned}
t_2 &= H(p, u) \\
&= p + (u-1)H_1 + (u-1)^2 H_2 + (p-P)(u-1)H_3 + (p-P)^2 H_4 + \ldots
\end{aligned}
$$ (2.6)

where,

$$H_1 = \frac{\partial H}{\partial u}\bigg|_{p=P, u=1}, \qquad\qquad H_2 = \frac{1}{2}\frac{\partial^2 H}{\partial u^2}\bigg|_{p=P, u=1},$$

$$H_3 = \frac{1}{2}\frac{\partial^2 H}{\partial p \partial u}\bigg|_{p=P, u=1}, \qquad \text{and} \qquad H_4 = \frac{1}{2}\frac{\partial^2 H}{\partial p \overline{y}^2}\bigg|_{p=P, u=1}.$$

The bias and MSE of the estimator $t_2$ are respectively given by –

$$B(t_2) = f\left(P\rho_{pb}C_p C_x H_3 + C_x^2 H_2 + P^2 C_y^2 H_4\right)$$ (2.7)

$$\text{MSE}(t_2) = f\left(P^2 C_p^2 + H_1^2 C_x^2 + 2H_1 P \rho_{pb} C_p C_x\right) \tag{2.8}$$

On differentiating (2.8) with respect to $H_1$ and equating to zero we obtain

$$H_1 = -\rho_{pb} P \frac{C_p}{C_x} \tag{2.9}$$

On substituting (2.9) in (2.8), we obtain the minimum MSE of the estimator $t_2$ as

$$\min \text{MSE}(t_2) = f P^2 C_p^2 \left(1 - \rho_{pb}^2\right) \tag{2.10}$$

We suggest another family of estimators for estimating P as

$$t_3 = \left[q_1 P + q_2 \left(\overline{X} - \overline{x}\right)\right]\left[\frac{a\overline{X} + b}{a\overline{x} + b}\right]^{\alpha} \exp\left[\frac{\left(a\overline{X} + b\right) - \left(a\overline{x} + b\right)}{\left(a\overline{X} + b\right) + \left(a\overline{x} + b\right)}\right]^{\beta} \tag{2.11}$$

where $\alpha, \beta, q_1$ and $q_2$ are real constants and a and b are known as characterising positive scalars. Many ratio-product estimators can be generated from $t_3$ by putting suitable values of $q_1$, $q_2$, $\alpha$, $\beta$, a and b (for choice of the parameters refer to [8], and [5]).

$$t_3 = \left[q_1 P(1 + e_0) - q_2 \overline{X}\right]\left[1 - \alpha\theta e_1 + \frac{\alpha(\alpha+1)}{2}\theta^2 e_1^2\right]\left[1 - \frac{\beta\theta e_1}{2} + \frac{\beta\theta^2 e_1^2}{8}(\beta+2)\right]$$

$$= q_1 P\left\{1 + e_0 - B(e_1 + e_0 e_1) + A e_1^2(1 + e_0)\right\} - q_2 \overline{X}\left\{e_1 - B e_1^2 \ldots\right\} \tag{2.12}$$

where, $\theta = \dfrac{a\overline{X}}{a\overline{X} + b}$, $B = \left(\alpha + \dfrac{\beta}{2}\right)\theta$ and $A = \dfrac{\theta^2}{8}\left[4\alpha(\alpha+1) + \beta(\beta+2) + 4\alpha\beta\right]$.

The bias and MSE of the estimator $t_3$ to the first order of approximation, are given as

$$\text{Bias}(t_3) = P(q-1) + f\left[\left(q_2 \overline{X}B + q_1 PA\right)C_x^2 - q_1 PB\rho C_p C_x\right] \tag{2.13}$$

$$MSE(t_3) = E(t_3 - P)^2$$

$$= (q_1 - 1)^2 P^2 + q_1^2 (M_1 + 2M_3) + q_2^2 M_2$$
$$+ 2q_1 q_2 (-M_4 - M_5) - 2q_1 M_3 + 2q_2 M_5 \qquad (2.14)$$

where,

$$M_1 = P^2 f\left(C_p^2 + B^2 C_x^2 - 2B\rho C_p C_x\right), \qquad M_2 = \overline{X}^2 f\left(C_x^2\right),$$

$$M_3 = P^2 f\left(AC_x^2 - 2B\rho C_p C_x\right), \qquad M_4 = P\overline{X} f\left(-BC_x^2 + \rho C_p C_x\right),$$

$$M_5 = \overline{X} P f\left(-BC_x^2\right)$$

On minimising the MSE of $t_3$ with respect to $q_1$ and $q_2$ respectively, we get

$$q_1^* = \frac{\Delta_1 \Delta_4 - \Delta_2 \Delta_5}{\Delta_1 \Delta_3 - \Delta_2^2} \qquad \text{and} \qquad q_2^* = \frac{\Delta_1 \Delta_5 - \Delta_2 \Delta_4}{\Delta_1 \Delta_3 - \Delta_2^2} \qquad (2.15)$$

where,

$$\Delta_1 = \left(P^2 + M_1 + 2M_3\right), \qquad \Delta_2 = (-M_4 - M_5),$$

$$\Delta_3 = (M_2), \qquad \Delta_4 = \left(P^2 + M_3\right),$$

$$\Delta_5 = (-M_5),$$

On putting these values of $q_1$ and $q_2$ in equation (2.14) we obtain the minimum MSE of $t_3$ as-

$$MSE(t_3)_{min} = \left[P^2 - \frac{\Delta_1 \Delta_5^2 + \Delta_3 \Delta_4^2 - 2\Delta_2 \Delta_4 \Delta_5}{\Delta_1 \Delta_3 - \Delta_2^2}\right] \qquad (2.16)$$

## 3. Efficiency Comparisons

First, we compare the efficiency of proposed estimator $t_3$ with usual estimator.

$$MSE(t_3)_{min} \leq V(\bar{y})$$

If,

$$\left[ P^2 - \frac{\Delta_1 \Delta_5^2 + \Delta_3 \Delta_4^2 - 2\Delta_2 \Delta_4 \Delta_5}{\Delta_1 \Delta_3 - \Delta_2^2} \right] \leq P^2 f_1 C_p^2 \qquad (3.1)$$

On solving we observed that above conditions holds always true.

Next we compare the efficiency of proposed estimator $t_3$ with regression estimator.

$$MSE(reg)MSE(t_\alpha)_{min} \leq MSE(reg)$$

If,

$$\left[ P^2 - \frac{\Delta_1 \Delta_5^2 + \Delta_3 \Delta_4^2 - 2\Delta_2 \Delta_4 \Delta_5}{\Delta_1 \Delta_3 - \Delta_2^2} \right] \leq P^2 f_1 C_p^2 \left( 1 - \rho_{pb}^2 \right) \qquad (3.2)$$

## 4. Empirical study

**Data Statistics:** We have taken the data from [1].

Where

Y – Home ownership

X – Income (thousands of dollars)

| n | N | P | $\overline{X}$ | $\rho_{pb}$ | Cp | Cx |
|---|---|---|---|---|---|---|
| 11 | 40 | 0.525 | 14.4 | 0..897 | 0.963 | 0.3085 |

The following Table shows PRE of different estimator's with respect to usual estimator.

**Table 1:   Percent relative efficiency (PRE) of estimators with respect to usual estimator**

| Estimators | $\overline{y}$ | $t_1$ | $t_2$ | $t_3$ | | |
|---|---|---|---|---|---|---|
| | | | | $\alpha=1, \beta=1$ | $\alpha=1, \beta=0$ | $\alpha=0, \beta=1$ |
| PRE | 100 | 189.384 | 511.794 | 515.798 | 517.950 | **518.052** |

When we examine Table 1, we observe that the proposed estimators $t_1$, $t_2$ and $t_3$ all performs better than the usual estimator $\overline{y}$. Also, the proposed estimator $t_3$ is the best among the estimators considered in the paper for the choice $\alpha=0, \beta=1$.